\documentclass[11pt]{article}
\usepackage{latexsym}
\usepackage[all]{xy}
\usepackage{amsmath}
\usepackage{amssymb}
\usepackage{amsfonts}
\usepackage{color}
\usepackage{theorem}

\newcommand{\Ql}{\mathbb{Q}_{\ell}}
\newcommand {\Map} {\mathbb{R}\mathbf{Map}}

\newcommand {\T} {\mathbb{T}}

\newcommand {\Spec} {\mathbf{Spec}}
\newcommand  {\dg}     {\mathbf{dg}}
\newcommand  {\dgcat}     {\mathbf{dgCat}}

\newcommand{\s}{\infty}

\newcommand{\BU}{\mathbb{B}\mathbb{U}}
\newcommand{\E}{\mathbb{E}}

\newcommand{\SH}{\mathcal{SH}}

\newtheorem{thm}{Theorem}[subsection]
\newtheorem{prop}[thm]{Proposition}

\newtheorem{df}[thm]{Definition}
\newtheorem{cor}[thm]{Corollary}
\newtheorem{conj}[thm]{Conjecture}

{\theorembodyfont{\rmfamily} \newtheorem{rmk}[thm]{Remark}}

\begin{document}

\title{\textbf{The $\ell$-adic trace formula for dg-categories \\ and Bloch's conductor conjecture}}  

\author{Bertrand To\"en\footnote{Partially supported by ANR-11-LABX-0040-CIMI within the program ANR-11-IDEX-0002-02 }  \, and Gabriele Vezzosi}  

\date{May 2016}

\maketitle

\begin{abstract} \noindent Building on \cite{brtv}, we present an $\ell$-adic trace formula for
smooth and proper dg-categories over a base $\E_\infty$-algebra $B$. We also give a variant when $B$ is only an $\E_2$-algebra.
As an application of this trace formula, we propose a strategy of proof of Bloch's conductor conjecture. \\
This is a research announcement and detailed proofs will appear elsewhere.

\end{abstract}

\tableofcontents

\section*{Introduction} 
This is a research announcement of results whose details and proofs will appear in forthcoming papers.\\ 

The purpose of this note is to give a trace formula in \'etale cohomology
in the non-commutative setting. Our non-commutative spaces will be modelled by 
dg-categories, and our trace formula reads as an equality between on one side, 
a certain virtual number of fixed points of an endomorphism, and on the other side
the trace of the induced endomorphism on the $\ell$-adic cohomology realization spaces
recently introduced in \cite{brtv} (and following the main idea of topological
realizations of dg-categories of \cite{blanc}). As an application of our
trace formula, we sketch a proof of the Bloch's conductor formula \cite{bl} 
under the hypothesis that the monodromy is unipotent. We also present a 
strategy for the general case. 

In order to describe our construction, we start by recalling some constructions and results from \cite{brtv}. In particular, for $A$ an 
excellent ring, $B$ an $\E_{\s}$ $A$-algebra, and $\ell$ a prime invertible on $A$, we describe the 
construction and basic properties of the lax symmetric monoidal $\ell$\emph{-adic realization 
functors} 
$$r_{\ell} : \dgcat_A \longrightarrow \Ql(\beta)-\mathbf{Mod} \, , \,\,\,  r^{B}_{\ell} : \dgcat_B \longrightarrow r_{\ell}(B)-\mathbf{Mod}$$ such that the induced square commutes up
to a natural equivalence
$$\xymatrix{
\dgcat_B \ar[r]^-{r_{\ell}^B} \ar[d] & r_{\ell}(B)-\mathbf{Mod} \ar[d] \\
\dgcat_A \ar[r]_-{r_{\ell}} & \Ql(\beta)-\mathbf{Mod}.}$$
Here $\Ql(\beta):= \oplus_{n\in \mathbb{Z}}(\mathbb{Q}_{\ell}[2](1))^{\otimes ^n}$ is viewed as an $\E_{\infty}$-ring 
object in the symmetric monoidal $\s$-category $\mathsf{L}(S_{\acute{e}t},\ell)$ of $\ell$-adic \'etale sheaves on $S= \Spec \, A$, and $\dgcat_A$ (resp. $\dgcat_B$) denotes the $\s$-category of dg-categories over $A$ (resp. over $B$, see Definition \ref{d2}).\\
Building on these results, we construct a \emph{Chern character} $$Ch: \mathsf{HK} \to |r_{\ell}(-)|$$ as a natural lax symmetric monoidal transformation  between $\s$-functors from $\mathbf{dgCat}_B$ to the $\s$-category $\mathbf{Sp}$ of spectra. Here $\mathsf{HK}$ is the non-connective homotopy invariant algebraic K-theory of dg-categories over $B$ (\cite{ro, brtv}), and, for any $T \in \mathbf{dgCat}_B$, $|r_{\ell}(T)|$ denotes the spectrum associated to the the hyper-cohomology of $r_{\ell}(T)$ via the (spectral) Eilenberg-Mac Lane construction.\\
With the Chern character at our disposal, we prove an $\ell$\emph{-adic trace formula} in the following setting. Let $T$ be a smooth, proper dg-category over $B$, and $f: T \to T$ and endomorphism of $T$, whose associated perfect $T^o \otimes_B T$-module is denoted by $\Gamma_f$. If $T$ satisfies a certain tensor-admissibility property  with respect to the $\ell$ -adic realization $r_{\ell}$ (called $r^{\otimes}_{\ell}$-admissibility in the text, see Definition \ref{d4}), to the effect that the natural coevaluation map exhibits $r_{\ell}(T^o)$ as the dual of $r_{\ell}(T)$ as an  $r_{\ell}(B)$-module, then there is a well defined \emph{trace} $\mathrm{Tr}(f_{| r_{\ell}(T)}) \in \pi_0(|r_{\ell}(B)|)$ of the endomorphism $f$. Our $\ell$-adic trace formula then states that $$\mathrm{Tr}(f_{| r_{\ell}(T)})= Ch([\mathsf{HH}(T/B,\Gamma_f)]),$$ where  $[\mathsf{HH}(T/B,\Gamma_f)]$ is the class in $\mathsf{HK}_0(B)$ of the Hochschild homology  of $\Gamma_f$ (which is indeed a perfect $B$-dg module).\\
With an eye to prospective applications, we also establish a version of this $\ell$-adic trace formula when the base $B$ is only supposed to be an $\E_2$-algebra (Section \ref{E2}).\\

We conclude this note by describing a strategy to prove Bloch's conductor conjecture (\cite{bl}, \cite[Conj. 6.2.1]{ks}), by using our $\ell$-adic trace formula together with the main result  of \cite{brtv}. Here $S= \Spec A$ is a henselian trait, with perfect residue field $k$, and fraction field $K$, and $X$ is a regular scheme, equipped with a proper and flat map to S, such that the generic fiber $X_K$ is smooth over $K$. We fix a uniformizer $\pi$ of $A$, and regard it as a global function on $X$. Our strategy is divided in two steps. In the first one we explain how to obtain a conductor formula, i.e. an expression of the Artin conductor of $X/S$ in terms of the Euler characteristic of the Hochschild homology of the category $\mathsf{MF}(X, \pi)$ of \emph{matrix factorizations} for the LG pair $(X, \pi)$. This conductor formula is interesting in itself, and seems to be new. The second step then consists in identifying our conductor formula with Bloch's.\\

\noindent \textbf{Acknowledgments.} We are grateful to our co-authors Anthony Blanc and Marco Robalo, for stimulating discussions that led to the joint work \cite{brtv}, and subsequently to the present paper. We also wish to thank Takeshi Saito for a very useful email exchange, and the Max-Planck-Institut f\"ur Mathematik in Bonn for providing a perfect scientific environment while the mathematics related to this paper was conceived.\\

\noindent \textbf{Notations.} Except where specified otherwise, we will write $S=\Spec\, A$ for $A$ an excellent commutative ring. We denote by $\mathbf{Sp}$ the $\s$-category of spectra, and by $\mathbb{S}\in \mathbf{Sp}$ the sphere spectrum. All functors (e.g. tensor products) will be implicitly derived. 

\section{$\ell$-adic realizations of nc-spaces}

\subsection{Reminders on dg-categories and their realizations}

We consider the category $dgCat_{A}$ of small $A$-linear dg-categories and $A$-linear dg-functors.
We remind that an $A$-linear dg-functor $T \longrightarrow T'$ is a Morita equivalence
if  the induced functor of the corresponding derived categories of dg-modules
$f^* : D(T') \longrightarrow D(T)$
is an equivalence of categories (see \cite{to-dgcat} for details). The $\s$-category of dg-categories over $S$
is defined to be the localisation of $dgCat_{A}$ along these Morita equivalences, 
and will be denoted by $\dgcat_S$ or $\dgcat_A$.

As in \cite[\S \, 4]{to-dgcat}, the tensor product of $A$-linear dg-categories can be derived to 
a symmetric monoidal structure on the $\s$-category $\dgcat_A$. It is a well known fact that 
dualizable objects in $\dgcat_S$ are precisely smooth and proper
dg-categories over $A$ (see \cite[Prop. 2.5]{to2}).

The compact objects in $\dgcat_S$ are the dg-categories of finite type over $A$ in the
sense of \cite{tv}. We denote their full sub-$\s$-category by $\dgcat_S^{ft} \subset \dgcat_S$.
The full sub-category $\dgcat_S^{ft}$ is preserved by  
the monoidal structure, and moreover any dg-category is a filtered colimit of 
dg-categories of finite type. 
We thus have a natural equivalence of symmetric monoidal $\s$-categories
$$\dgcat_S \simeq \mathbf{Ind}(\dgcat_S^{ft}).$$

We denote by $\SH_S$ the stable $\mathbb{A}^1$-homotopy $\s$-category of schemes over $S$ (see \cite[Def. 5.7]{vo} and \cite[\S\, 2]{ro}). 
It is a  presentable symmetric monoidal $\s$-category whose monoidal structure will be denoted by 
$\wedge_S$. Homotopy invariant algebraic K-theory defines an $\E_{\infty}$-ring 
object  in $\SH_S$ that we denote by $\BU_S$ (a more standard notation is $KGL$). We denote by 
$\BU_S-\mathbf{Mod}$ the $\s$-category of $\BU_S$-modules in $\SH_S$. It is a presentable symmetric monoidal
$\s$-category whose monoidal structure will be denoted by $\wedge_{\BU_S}$. 

As proved in \cite{brtv}, there exists a lax symmetric monoidal $\s$-functor
$$\mathsf{M}^{-}: \dgcat_S \longrightarrow \BU_S-\mathbf{Mod},$$
which is denoted by $T \mapsto \mathsf{M}^T$. The precise construction of the $\s$-functor $\mathsf{M}^{-}$
is rather involved and uses in an essential manner the theory of non-commutative motives 
of \cite{ro} as well as the comparison with the stable homotopy theory of schemes. Intuitively, 
the $\s$-functor $\mathsf{M}^{-}$ sends a dg-category 
 $T$ to the homotopy invariant K-theory functor $S' \mapsto \mathsf{HK}(S' \otimes_S T)$. 
To be more precise, there is an obvious forgetful $\s$-functor 
$$\mathsf{U}: \BU_S-\mathbf{Mod} \longrightarrow \mathbf{Fun}^{\s}(Sm_S^{op},\mathbf{Sp}),$$
to the $\s$-category of presheaves of spectra on the category $Sm_S$ of smooth $S$-schemes. For a given 
dg-category $T$ over $S$, the presheaf $\mathsf{U}(\mathsf{M}^T)$ is defined 
 by sending a smooth $S$-scheme $S'=\Spec\, A' \rightarrow \Spec\, A=S$
to $\mathsf{HK}(A'\otimes_A T)$, the homotopy invariant non-connective K-theory spetrum
of $A'\otimes_A T$ (see \cite[4.2.3]{ro}). \\

The $\s$-functor $\mathsf{M}^{-}$ satisfies some basic properties which we recall here.
\begin{enumerate}
\item The $\s$-functor $\mathsf{M}^{-}$ is a localizing invariant, i.e. for any 
short exact sequence $T_0 \hookrightarrow T \longrightarrow T/T_0$ of dg-categories over $A$, the induced sequence 
$$\xymatrix{\mathsf{M}^{T_0} \ar[r] & \mathsf{M}^T \ar[r] & \mathsf{M}^{T/T_0}}$$
exhibits $\mathsf{M}^{T_0}$ has the fiber of the morphism $\mathsf{M}^T \rightarrow \mathsf{M}^{T/T_0}$ in $\BU_S-\mathbf{Mod}$.

\item The natural morphism
$\BU_S \longrightarrow \mathsf{M}^{A},$
induced by the lax monoidal structure of $\mathsf{M}^{-}$, is an equivalence of $\BU_S$-modules. 

\item The $\s$-functor $T \mapsto \mathsf{M}^T$ commutes with filtered colimits.

\item For any quasi-compact and quasi-separated scheme $X$, and any morphism  $p : X \longrightarrow S$, 
we have a natural equivalence of $\BU_S$-modules
$$\mathsf{M}^{\mathsf{L}_{\textrm{Perf}}(X)} \simeq p_*(\BU_X),$$
where $p_* : \BU_X-\mathbf{Mod} \longrightarrow \BU_S-\mathbf{Mod}$ is the direct image of
$\BU$-modules, and $\mathsf{L}_{\textrm{Perf}}(X)$ is the dg-category of perfect complexes on $X$.

\end{enumerate}

We now let $\ell$ be a prime number invertible in $A$. We denote by 
$\mathsf{L}_{\mathsf{ct}}(S_{\acute{e}t},\ell)$ the $\s$-category of constructible 
$\mathbb{Q}_{\ell}$-complexes on the \'etale site $S_{\acute{e}t}$ of $S$. It is a symmetric monoidal $\s$-category, and
we denote by 
$$\mathsf{L}(S_{\acute{e}t},\ell):=\mathbf{Ind}(\mathsf{L}_{\mathsf{ct}}(S_{\acute{e}t},\ell))$$
its completion under filtered colimits (see \cite[Def. 4.3.26]{gl}). 
According to \cite[Cor. 2.3.9]{ro}, there exists an $\ell$-adic realization $\s$-functor
$r_{\ell} : \SH_S \longrightarrow \mathsf{L}(S_{\acute{e}t},\ell).$
By construction, $r_{\ell}$ is a symmetric monoidal $\s$-functor 
sending a smooth scheme $p : X \longrightarrow S$ to $p_!p^{!}(\mathbb{Q}_{\ell})$, 
or, in other words, to the relative $\ell$-adic homology of $X$ over $S$. 

We let $\T:=\mathbb{Q}_{\ell}[2](1)$, and we consider the $\E_{\infty}$-ring 
object in $\mathsf{L}(S_{\acute{e}t},\ell)$
$$\Ql(\beta):=\oplus_{n\in \mathbb{Z}}\T^{\otimes n}.$$
 In this notation, $\beta$ stands for $\T$, 
and $\Ql(\beta)$ for the algebra of Laurent polynomials in $\beta$, so we could
also have written
$$\Ql(\beta)=\Ql[\beta,\beta^{-1}].$$
As shown in \cite{brtv}, there exists a canonical equivalence $r_{\ell}(\BU_S) \simeq \Ql(\beta)$ of $\E_{\s}$-ring objects
in $\mathsf{L}(S_{\acute{e}t},\ell)$,
that is induced by the Chern character from algebraic K-theory to 
\'etale cohomology. We thus obtain a well-defined symmetric monoidal $\s$-functor
$$r_{\ell} : \BU_S-\mathbf{Mod} \longrightarrow \Ql(\beta)-\mathbf{Mod},$$
from $\BU_S$-modules in $\SH_S$ to $\Ql(\beta)$-modules in $\mathsf{L}(S_{\acute{e}t},\ell)$.
By pre-composing with the functor $T \mapsto \mathsf{M}^T$, we obtain
a lax monoidal $\s$-functor 
$$r_{\ell} : \dgcat_S \longrightarrow \Ql(\beta)-\mathbf{Mod}.$$

\begin{df}\label{d1}
The $\s$-functor defined above 
$$r_{\ell} : \dgcat_S \longrightarrow \Ql(\beta)-\mathbf{Mod}$$
is called the \emph{$\ell$-adic realization functor for dg-categories over $S$}.
\end{df}  

From the standard properties of the functor $T \mapsto \mathsf{M}^T$, recalled above, we obtain the following
properties for the $\ell$-adic realization functor $T \mapsto r_{\ell}(T)$. 
 
\begin{enumerate}
\item The $\s$-functor $r_{\ell}$ is a localizing invariant, i.e. for any 
short exact sequence $T_0 \hookrightarrow T \longrightarrow T/T_0$ of dg-categories over $A$, the induced sequence 
$$\xymatrix{r_{\ell}(T_0) \ar[r] & r_{\ell}(T) \ar[r] & r_{\ell}(T/T_0)}$$
is a fibration sequence in $\Ql(\beta)-\mathbf{Mod}$.

\item The natural morphism
$$\Ql(\beta) \longrightarrow r_{\ell}(A),$$
induced by the lax monoidal structure, is an equivalence in $\Ql(\beta)-\mathbf{Mod}$. 

\item The $\s$-functor $r_{\ell}$ commutes with filtered colimits.

\item For any separated morphism of finite type  $p : X \longrightarrow S$, 
we have a natural morphism of $\Ql(\beta)$-modules
$$r_{\ell}(\mathsf{L}_{\textrm{Perf}}(X)) \longrightarrow p_*(\Ql(\beta)),$$
where $p_* : \Ql(\beta)-\mathbf{Mod} \longrightarrow \Ql(\beta)-\mathbf{Mod}$ is induced by the
direct image $L_{\mathsf{ct}}(X_{\acute{e}t},\ell) \longrightarrow L_{\mathsf{ct}}(S_{\acute{e}t},\ell)$ of 
constructible $\Ql$-complexes.
If $p$ is proper, or if $A$ is a field, this morphism is an equivalence.

\end{enumerate}
 
To finish this part, we generalize a bit the above setting, by adding to the context a base
$\E_\infty$-algebra $B$ over $A$, and considering $B$-linear dg-categories instead of just $A$-linear dg-categories.\\

Let $B$ be an $\E_\infty$-algebra over $A$. We consider $B$ as
an $\E_\infty$-monoid in the symmetric monoidal $\s$-category $\dgcat_A$.  
We define dg-categories over $B$ as being $B$-modules in $\dgcat_A$. More
specifically we have the following notion.

\begin{df}\label{d2}
The \emph{symmetric monoidal $\s$-category of (small) $B$-linear dg-categories} is defined to be the $\s$-category
of $B$-modules in $\dgcat_A$. It is denoted by
$$\dgcat_B:=B-\mathbf{Mod}_{\dgcat_A},$$
and its monoidal structure by $\otimes_B$.
\end{df}
 
By applying our $\ell$-adic realization functor (Definition \ref{d1}), we have that $r_{\ell}(B)$ is
an $\E_\infty$ $\Ql(\beta)$-algebra. We thus get an induced lax symmetric monoidal $\s$-functor
$$r^{B}_{\ell} : \dgcat_B \longrightarrow r_{\ell}(B)-\mathbf{Mod}.$$
By construction, the natural forgetful $\s$-functors make the following square commute up
to a natural equivalence
$$\xymatrix{
\dgcat_B \ar[r]^-{r_{\ell}^B} \ar[d] & r_{\ell}(B)-\mathbf{Mod} \ar[d] \\
\dgcat_A \ar[r]_-{r_{\ell}} & \Ql(\beta)-\mathbf{Mod}.}$$

We will often consider $B$ as implicitly assigned, and we will simply write $r_{\ell}$ for
$r_{\ell}^B$.
  
\subsection{Chern character and Grothendieck-Riemann-Roch}

We fix as above an $\E_\infty$-algebra $B$ over $A$. As explained in 
\cite{brtv}, there is a symmetric monoidal $\s$-category $\SH_B^{nc}$ of non-commutative motives over $B$. 
As an $\s$-category it is the full sub-$\s$-category of $\s$-functors
of (co)presheaves of spectra
$$\dgcat_B^{ft} \longrightarrow \mathbf{Sp},$$
satisfying Nisnevich descent and $\mathbb{A}^1$-homotopy invariance. The symmetric monoidal
structure is induced by left Kan extension of the symmetric monoidal
structure $\otimes_B$ on $\dgcat_B^{ft}$. 

We consider 
$\Gamma : \mathsf{L}(S_{\acute{e}t},\ell) \longrightarrow \dg_{\Ql},$
the global section $\s$-functor, taking 
an $\ell$-adic complex on $S_{\acute{e}t}$ to its hyper-cohomology. Composing this 
with the Dold-Kan construction $\Map_{\dg_{\Ql}}(\Ql,-) : \dg_{\Ql} \longrightarrow \mathbf{Sp}$, 
we obtain an $\s$-functor
$$|-| : \mathsf{L}(S_{\acute{e}t},\ell) \longrightarrow \mathbf{Sp},$$
which morally computes hyper-cohomology of $S_{\acute{e}t}$ with $\ell$-adic coefficients, i.e. 
for any $E \in \mathsf{L}(S_{\acute{e}t},\ell)$, we have natural isomorphisms
$$H^i(S_{\acute{e}t},E) \simeq \pi_{-i}(|E|)\,, i\in \mathbb{Z}.$$
By what we have seen in our last paragraph, the composite functor 
$T \mapsto |r_{\ell}(T)|$ provides a
(co)presheaves of spectra
$$\dgcat_B^{ft} \longrightarrow \mathbf{Sp},$$
satisfying Nisnevich descent and $\mathbb{A}^1$-homotopy invariance.
It thus defines an object $|r_{\ell}| \in \SH_B^{nc}$. The fact that $r_{\ell}$
is lax symmetric monoidal implies moreover that $|r_{\ell}|$ is endowed with 
a natural structure of a $\E_\infty$-ring object in $\SH_B^{nc}$. 

Each $T \in \dgcat_B^{ft}$ defines a corepresentable
object $h^T \in \SH_B^{nc}$, characterized by the ($\s$-)functorial equivalences
$$\Map_{\SH_B^{nc}}(h^T,F) \simeq F(T),$$
for any $F \in \SH_B^{nc}$. The existence of $h^T$ is a formal statement, however
the main theorem of \cite{ro} implies that we have natural equivalences of spectra
$$\Map_{\SH_B^{nc}}(h^T,h^B) \simeq \mathsf{HK}(T),$$
where $\mathsf{HK}(T)$ stands for non-connective homotopy invariant algebraic $K$-theory of the dg-category $T$.
In other words, $T \mapsto \mathsf{HK}(T)$ defines an object in $\SH_B^{nc}$ which 
is isomorphic to $h^B$.
By Yoneda lemma, we thus obtain an equivalence of spaces
$$\Map^{lax-\otimes}(\mathsf{HK},|r_\ell|) \simeq \Map_{\E_\infty-\mathbf{Sp}}(\mathbb{S},|r_{\ell}(B)|)\simeq *.$$
In other words, there exists a unique (up to a contractible space of choices) lax symmetric monoidal natural transformation 
$$\mathsf{HK} \longrightarrow |r_\ell|,$$
between lax monoidal $\s$-functors from $\dgcat^{ft}_B$ to $\mathbf{Sp}$. We 
extend this to all dg-categories over $B$ as usual by passing to Ind-completion
$\dgcat_B \simeq \mathbf{Ind}(\dgcat_B^{ft})$.

\begin{df}\label{d3}
The natural transformation defined above is called the \emph{$\ell$-adic Chern character}.
It is denoted by 
$$Ch : \mathsf{HK}(-) \longrightarrow |r_{\ell}(-)|.$$ 
\end{df}

Note that the Chern character is an absolute construction, and does not depend on the base
$B$. In other words, it factors through the natural forgetful functor $\dgcat_B \longrightarrow \dgcat_A$.
There is thus no need to specify the base $B$ in the notation $Ch$. \\

Definition \ref{d3} contains a formal Grothendieck-Riemann-Roch formula. Indeed, 
for any $B$-linear dg-functor $f : T \longrightarrow T'$, the square of spectra
$$\xymatrix{
\mathsf{HK}(T) \ar[r]^-{f_!} \ar[d]_-{Ch_T} & \mathsf{HK}(T') \ar[d]^-{Ch_{T'}} \\
r_{\ell}(T) \ar[r]_-{f_!} & r_{\ell}(T')}$$
commutes up to a natural equivalence.

\subsection{The trace formula}

In the last paragraph of this section we prove an $\ell$\emph{-adic trace formula} for
smooth and proper dg-categories over a base $\E_\infty$-algebra $B$. The formula 
is a direct consequence of dualizability of smooth and proper dg-categories, 
and functoriality of the Chern character. 

We let $B$ be a fixed $\E_\infty$-algebra over $A$ and let
$\dgcat_B$ the symmetric monoidal $\s$-category of dg-categories over $B$. Recall
that the dualizable objects $T$ in $\dgcat_B$ are precisely  
smooth and proper dg-categories (see \cite[Prop. 2.5]{to2}). For such dg-categories $T$, the dual $T^{\vee}$ in $\dgcat_B$ is exactly the opposite dg-category $T^o$, and the evaluation map
$$T^o \otimes_B T \longrightarrow B$$
is simply given by sending $(x,y)$ to the perfect $B$-module $T(x,y)$. In the
same manner, the coevaluation map
$$B \longrightarrow T^o \otimes_B T$$
is given the identity bi-module $T$. In particular, the
composition
$$\xymatrix{B \ar[r] & T^o\otimes_B T \ar[r] & B}$$
is the endomorphism of $B$ given by the perfect $B$-dg module $\mathsf{HH}(T/B)$, the Hochschild
chain complex of $T$ over $B$. 

\begin{df}\label{d4}
We say that a smooth and proper dg-category $T$ over $B$ is 
\emph{$r_{\ell}^{\otimes}$-admissible} if  the natural morphism, induced by the
lax monoidal structure
$$r_{\ell}(T^o)\otimes_{r_{\ell}(B)}r_{\ell}(T) \longrightarrow r_\ell(T^o \otimes_B T)$$
is an equivalence in $\mathsf{L}(S_{\acute{e}t},\ell)$.
\end{df}

A direct observation shows that if $T \in \dgcat_B$ is smooth, proper and \emph{$r_{\ell}^{\otimes}$-admissible}, 
then the image by $r_{\ell}$ of the co-evaluation map
$$r_{\ell}(B) \longrightarrow r_{\ell}(T^o \otimes_B T) \simeq r_{\ell}(T^o) \otimes_{r_{\ell}(B)}r_{\ell}(T)$$
exhibits $r_{\ell}(T^o)$ as the dual of $r_{\ell}(T)$ as a $r_{\ell}(B)$-module. \\

\noindent We now fix a smooth, proper and $r_{\ell}^{\otimes}$-admissible dg-category $T$ over $B$. Let $f : T \longrightarrow T$
be an endomorphism of $T$ in $\dgcat_B$. By duality, $f$ corresponds to 
a perfect $T^o\otimes_{B}T$-dg module denoted by 
$\Gamma_f$. The Hochschild homology of $\Gamma_f$ is then the perfect $B$-dg module
$\mathsf{HH}(T/B,\Gamma_f)$, defined as the composite morphism in $\dgcat_B$
$$\xymatrix{
B \ar[r]^-{\Gamma_f} & T^o \otimes_B T \ar[r]^-{\mathsf{ev}} & B.}$$
Since $\mathsf{HH}(T/B,\Gamma_f)$ is a  perfect $B$-dg module, it defines a class  
$$[\mathsf{HH}(T/B,\Gamma_f)] \in \mathsf{HK}_0(B).$$

On the other hand, the given endomorphism $f$ induces by functoriality an endomorphism of $r_{\ell}(B)$-modules
$$f_{|r_{\ell}(T)}:= r_{\ell}^B(f) : r^B_{\ell}(T) \longrightarrow r^B_{\ell}(T).$$
Since we are supposing that $T$ is $r_{\ell}^{\otimes}$-admissible, as already observed $r^B_{\ell}(T)$ is dualizable (with dual $r^B_{\ell}(T^o)$), and therefore we are able to define the trace of $f_{|r_{\ell}(T)}$, again, as the composite
$$\xymatrix{
r_{\ell}(B) \ar[rr]^-{\Gamma_{f_{|r_{\ell}(T)}}} & & r^B_{\ell}(T)^{\vee} \otimes_{r_{\ell}(B)}r^B_{\ell}(T) \ar[r]^-{\mathsf{ev}} & 
r_{\ell}(B).}$$
The equivalence class of this morphism is identified with an element
$$\mathrm{Tr}(f_{|r_{\ell}(T)}) \in \pi_0(|r_{\ell}(B)|).$$

Our $\ell$-adic trace formula, which is a direct consequence of duality in $\dgcat_B$ and in $r_{\ell}(B)-\mathbf{Mod}$, is the following statement

\begin{thm}\label{t1}
Let $T \in \dgcat_B$  be a smooth, proper and $r_{\ell}^{\otimes}$-admissible dg-category over $B$, and $f : T \longrightarrow T$ be an endomorphism. Then we have
$$Ch([\mathsf{HH}(T/B,\Gamma_f)])=\mathrm{Tr}(f_{|r_{\ell}(T)}).$$
\end{thm}

A case of special interests is when $B$ is such that the natural morphisms of rings
$$\mathbb{Z} \longrightarrow \mathsf{HK}_0(B) \qquad \Ql \longrightarrow \pi_0(|r_{\ell}(B)|)$$
are both isomorphisms. In this case the Chern character
$Ch : \mathsf{HK}_0(B) \longrightarrow H^0(|r_\ell(B)|)$
is the natural inclusion $\mathbb{Z} \subset \Ql$, and the formula reads just 
as an equality of $\ell$-adic numbers
$$[\mathsf{HH}(T/B,\Gamma_f)]=\mathrm{Tr}(f_{|r_{\ell}(T)}).$$

The left hand side of this formula should be interpreted as the intersection number
of the graph $\Gamma_f$ with the diagonal of $T$, and thus
as a \emph{virtual number of fixed points of} $f$.

\subsection{Extension to $\E_2$-bases}\label{E2}

The trace formula of Theorem \ref{t1} can be extended to the case
where the base $\E_\infty$-algebra $B$ is only 
assumed to be an $\E_2$-algebra. We will here briefly sketch how this works. \\

Let $B$ be an $\E_2$-algebra over $A$. It can be considered 
as an $\E_1$-algebra object in the symmetric monoidal $\s$-category $\dgcat_A$, and  
we thus define the $\s$-category of dg-categories enriched over $B$ (or linear over $B$) as the $\s$-category 
$$\dgcat_B:=B-\mathbf{Mod}_{\dgcat_A},$$
 of $B$-modules in $\dgcat_A$. As opposed to the case where $B$ is $\E_\infty$, 
$\dgcat_B$ is no more a symmetric monoidal $\s$-category, and dualizability must 
therefore be understood in a slightly different manner.

We will start by working in $\dgcat^{big}_A$, an enlargement of 
$\dgcat_A$ where we allow for arbitrary bimodules as morphisms. More precisely, 
$\dgcat^{big}_A$ is an $\s$-category whose objects are small dg-categories
over $A$, and morphisms between $T$ and $T'$ in $\dgcat_A^{big}$ are
given by the space of all $T\otimes_A (T')^{o}$-dg modules. Equivalently, 
$\dgcat_A^{big}$ is the $\s$-category of all compactly generated
presentable dg-categories over $A$ and continuous dg-functors. 
The tensor product of dg-categories endows
$\dgcat_A^{big}$ with a symmetric monoidal structure 
still denoted by $\otimes_A$.
We have a natural symmetric monoidal $\s$-functor $$\dgcat_A \longrightarrow \dgcat_A^{big},$$
identifying $\dgcat_A$ with the sub-$\s$-category of 
continuous morphisms preserving compact objects. As these two 
$\s$-categories only differ by their morphisms, we will write \emph{big morphisms} to mean morphisms in $\dgcat^{big}_A$. \\

Going back to our $\E_2$-algebra $B$, the product $m_B: B\otimes_A B^o \longrightarrow B$ is a morphism of $\E_1$-algebras, and we may consider the composite functor 
$$\xymatrix{\delta_B: B\otimes_A B^o -\mathbf{mod}\ar[r]^-{m^*_B} & B-\mathbf{mod} \ar[r]^-{u_*} & A-\mathbf{mod}}$$ between the corresponding categories of \emph{dg-modules} (not of dg-categories), where $u_*$ denotes the forgetful functor induced by the canonical map $u: A \to B$. Note that, $B$ being an $\E_2$-algebra, the composite functor $\delta_B$ is lax monoidal (even though the base-change $m^{*}_B$, alone, is not), and thus induces a well defined functor on the corresponding $\s$-categories of modules
in $\dgcat_A^{big}$, denoted by 
$$\Delta_{B} : B\otimes_A B^o-\mathbf{Mod} \longrightarrow A-\mathbf{Mod}=\dgcat_A^{big},$$ where we simply write $\mathbf{Mod}$ for $\mathbf{Mod}_{\dgcat_A^{big}}$. Note that, by definition, $$\Delta_{B}(B\otimes_A B^o)=B,$$ where $B$ viewed as an object in $\dgcat_A^{big}$.
As a consequence, for any $T \in \dgcat_B$, we may define $T\otimes_B T^o$ by the formula
$$T\otimes_B T^o := \Delta_B(T\otimes_A T^o) \in \dgcat_A.$$
The $B$-module structure morphism on $T$
$$B\otimes_A T \longrightarrow T$$
provides a big morphism
$$\mathsf{ev} : T\otimes_A T^o \longrightarrow B^o,$$
which is a morphism of $B\otimes_A B^o$-modules. By applying $\Delta_B$, we get an induced \emph{evaluation} morphism in $\dgcat^{big}_A$
$$\mathsf{ev}' :  T\otimes_B T^o \longrightarrow \Delta_{B}(B^o)=:\mathsf{HH}(B/A).$$
Notice that $\mathsf{ev}'$ is only a big morphism of $A$-linear dg-categories as there is no natural 
$B$-linear structure on $T\otimes_B T^o$. Moreover, since $B$ is an $\E_2$-algebra,  $\mathsf{HH}(B/A)\simeq S^1 \otimes_A B$ is itself 
an $\E_1$-algebra and thus can be considered as an object of $\dgcat_A$.\\
We leave to the reader the, similar, dual construction of a \emph{coevaluation} big morphism
$$\mathsf{coev} : A \longrightarrow T^o\otimes_B T.$$
More generally, for $f : T \longrightarrow T$ an endomorphism of $T$ in $\dgcat_B$,
we have an induced graph morphism in $\dgcat_A^{big}$
$$\Gamma_f : A \longrightarrow T^o \otimes_B T.$$

\begin{df}\label{d5}
\begin{enumerate}
\item
A $B$-linear dg-category $T\in \dgcat_B$ is \emph{smooth (resp. proper)} if 
the coevaluation $\mathsf{coev}$ (resp. evaluation $\mathsf{ev}$) big morphism actually lies in $\dgcat_A$ (i.e. 
preserves compact objects).

\item For a smooth and proper $B$-linear dg-category $T$, and an
endomorphism $f : T \rightarrow T$ in $\dgcat_B$, the \emph{$\E_2$-Hochschild 
homology of $T$ relative to $B$ with coefficients in $f$} is defined 
to be the composite morphism in $\dgcat_A$
$$\mathsf{HH}^{\E_2}(T/B,f):=
\xymatrix{A \ar[r]^-{\mathsf{coev}} & T\otimes_B T^o \ar[r]^-{\mathsf{ev}'} & \mathsf{HH}(B/A).}$$
\end{enumerate}
\end{df}

Note that, by definition, the Hochschild homology of the pair $(T,f)$ 
can be identified with a perfect $\mathsf{HH}(B/A)$-module. As
$B$ is an $\E_2$-algebra, the natural augmentation
$\alpha: \mathsf{HH}(B/A) \longrightarrow B$ is an $\E_1$-morphism. Therefore, we can 
base-change $\mathsf{HH}(T/B,f)$ along $\alpha$, and get a perfect $B$-dgmodule
$\mathsf{HH}^{\E_2}(T/B,f) \otimes_{\mathsf{HH}(B/A)}B$. When $B$ is an $\E_\infty$-algebra
the Hochschild homology and the $\E_2$-Hochschild homology are related by 
a natural equivalence of $B$-dg-modules
$$\mathsf{HH}(T/B,f) \simeq \mathsf{HH}^{\E_2}(T/B,f)\otimes_{\mathsf{HH}(B/A)}B.$$
To be more precise, we have an equivalence 
$$\mathsf{HH}^{\E_2}(T/B,f) \simeq \mathsf{HH}(T/A,f)$$ of perfect $\mathsf{HH}(B/A)$-dg-modules.\\

For a smooth and proper $B$-linear dg-category $T$, we consider the
following (compare with Definition \ref{d4})

\begin{itemize} \item \textbf{$r_{\ell}^{\otimes}$-admissibility assumptions.} The following natural morphisms 
\begin{align*}  
{r_{\ell}(T)\otimes_{r_{\ell}(B)}r_{\ell}(T^o) \longrightarrow r_{\ell}(T\otimes_B T^o)\, , \qquad
\mathsf{HH}(r_{\ell}(B)/\Ql(\beta)) \longrightarrow r_{\ell}(\mathsf{HH}(B/A))} 
\tag{Adm}\end{align*} are equivalences 
in $\mathsf{L}(S_et,\ell)$.
\end{itemize}
Under these assumptions, the evaluation and coevaluation morphisms for $T$ induce well defined morphisms
$$\mathsf{coev} : \Ql(\beta) \longrightarrow
r_{\ell}(T)\otimes_{r_{\ell}(B)}r_{\ell}(T^o) \qquad 
\mathsf{ev} : r_{\ell}(T)\otimes_{\Ql(\beta)}r_{\ell}(T^o) \longrightarrow r_{\ell}(B).
$$
Here, the morphism $\mathsf{coev}$ is $\Ql(\beta)$-linear while $\mathsf{ev}$ is $r_{\ell}(B)\otimes_{\Ql(\beta)}r_{\ell}(B^o)$-linear.
These two morphisms exhibit $r_{\ell}(T^o)$ has the right dual of $r_{\ell}(T)$ in $r_{\ell}(B)$-modules. 
In particular, an endomorphism $f:T \to T$ in $\dgcat_B$  together with the evaluation morphism provide a well defined trace morphism
$$\mathrm{Tr}(f) : \Ql(\beta) \longrightarrow r_{\ell}(B)\otimes_{r_{\ell}(B)\otimes_{\Ql(\beta)}r_{\ell}(B^o)}r_{\ell}(B^o),$$
and a corresponding well defined element
$$\mathrm{Tr}(f_{|r_{\ell}(T)}) \in H^0(S_{\acute{e}t},\mathsf{HH}(r_{\ell}(B)/\Ql(\beta))).$$
The $\E_2$-version of the trace formula of Theorem \ref{t1} is then the following

\begin{thm}\label{t2}
Let $B$ be an $\E_2$-algebra over $A$, $T$ a smooth and proper $B$-linear dg category, and $f: T \to T$ a $B$-linear endomorphism. We have an equality 
$$Ch([\mathsf{HH}^{\E_2}(T/B,f)])=\mathrm{Tr}(f_{|r_{\ell}(T)})$$ in $H^0(S_{\acute{e}t},\mathsf{HH}(r_{\ell}(B)/\Ql(\beta)))$.
\end{thm}

By using the augmentation $\rho : \mathsf{HH}(B/A) \longrightarrow B$, we deduce a perhaps more familiar version of the trace formula, as 
$$Ch([\mathsf{HH}^{\E_2}(T/B,f)\otimes_{\mathsf{HH}(B/A)}B]) = \rho(\mathrm{Tr}(f_{|r_{\ell}(T)})),$$
which is now an equality in $H^0(S_{\acute{e}t},r_{\ell}(B))$.

\section{Applications to vanishing cycles and Bloch's conductor conjecture}

Our trace formula (Theorem \ref{t1} and Theorem \ref{t2}), combined with the main result of \cite{brtv}, has an interesting application
to Bloch's conductor conjecture of \cite{bl}. 

\subsection{Bloch's conductor conjecture} Our base scheme 
is an henselian trait $S=\Spec\, A$, with perfect residue field $k$, and fraction field $K$. 
Let $X \longrightarrow S$
be proper and flat morphism of finite type, and of relative dimension $n$. We assume that the generic fiber $X_K$ is smooth over $K$, and 
that $X$ is a regular scheme. Bloch's conductor formula conjecture
reads as follows (see \cite{bl,ks}, for detailed definitions of the various objects involved).

\begin{conj}\label{cb} \emph{\textbf{[Bloch's conductor Conjecture]}}
We have an equality
$$[\Delta_X,\Delta_X]_S = \chi(X_{\bar{k}},\ell) - \chi(X_{\bar{K}},\ell) - Sw(X_{\bar{K}}),$$
where $\chi(Y,\ell)$ denotes the $\ell$-adic Euler characteristic of a variety
$Y$ for $\ell$ prime to the characteristic of $k$, $Sw(X_{\bar{K}})$
is the Swan conductor of the $Gal(K)$-representation $H^*(X_{\bar{K}},\Ql)$,
and $[\Delta_X,\Delta_X]_S$ is the degree in $\textrm{CH}_0(k)\simeq \mathbb{Z}$ of Bloch's localised self-intersection $(\Delta_X,\Delta_X)_S \in \mathrm{CH}_{0}(X_k)$ of the
diagonal in $X$. The (negative of the) rhs is called the Artin conductor of $X/S$.
\end{conj}

Note that this conjecture is known in some important cases: it is classical for $n=0$, proved by Bloch himself for $n=1$ \cite{bl}, and for arbitrary $n$ by Kato and Saito \cite{ks} under the hypothesis that the reduced special fiber $(X_k)_{\textrm{red}}$ is a normal crossing divisor. By \cite{or}, we also know that Conjecture \ref{cb} implies the Deligne-Milnor Conjecture for isolated singularities (both in the equi and mixed characteristic case). The isolated singularities and equicharacteristic case of the Deligne-Milnor conjecture was proved by Deligne \cite[Exp. XVI]{sga7II}, while the isolated singularities and mixed characteristic case was proved by Orgogozo \cite{or} for relative dimension $n=1$. \\

\subsection{A strategy}
We propose here an approach to Bloch's conjecture based on our trace formula (Theorem \ref{t2}) and the use of the theory of \emph{matrix factorizations}, as developed in \cite{brtv}. 
Below we will just sketch a general strategy that might lead to a proof of the general case of Conjecture \ref{cb}. Some details are still to be fixed, namely the equality (F-gen) and the comparison formula (Comp) with Bloch's localised self-intersection of the diagonal (see below). These details are currently being checked and, hopefully, they will appear in a forthcoming paper. \\

Let $\pi$ be a uniformizer for $A$, so that $k=A/\pi$. We let $T:=\mathsf{MF}(X,\pi)$ be the $A$-linear dg-category of matrix factorizations
of $X$ for the function $\pi$ on $X$, as studied in detail in \cite{brtv}. We consider
the $\E_2$-algebra $B^{+}$ defined as 
$$B^{+}:=\mathbb{R}\underline{Hom}_{k\otimes_A k}(k,k).$$
The fact that $B^{+}$
is an $\E_2$-algebra follows formally from the fact that 
it is endowed with two compatible $\E_1$-multiplications, 
one given by composition of endomorphisms, and another
coming from the fact that the derived scheme $\Spec\, (k\otimes^{L}_A k)$ is actually a Segal groupoid
object in derived schemes: the descent 
groupoid of the map $\Spec\, k \rightarrow S$. 
The composition law in this groupoid induces another $\E_1$-structure on $B^{+}$ given 
by convolution, making it into an $\E_2$-algebra over $A$.

As an $\E_1$-algebra, $B^{+}$ can be identified with $k[u]$, where
$u$ is a free variable in degree $2$. We set
$$B:=B^{+}[u^{-1}].$$
Again, as an $\E_1$-algebra $B$ is just $k[u,u^{-1}]$. However, in general,
$B$ is not equivalent to $k[u,u^{-1}]$ as an $\E_2$-algebra, as, a priori, it is not even
linear over $k$. It can be shown that if $A$ is a $k$-algebra (i.e. we are in the geometric case), then 
$B$ is equivalent to $k[u,u^{-1}]$ and thus
is indeed an $\E_\infty$-algebra over $k$. \\

In \cite{brtv} it is shown that the $\E_2$-algebra $B$ naturally acts on the dg-category 
$T=\mathsf{MF}(X,\pi)$ making it into a $B$-linear dg-category (as in \S \ref{E2}). Note 
that we have a natural ring isomorphism 
$$\chi : \mathsf{HK}_0(B) \simeq \mathbb{Z}.$$

A key result in our strategy is the the following proposition, whose proof will appear elsewhere.

\begin{prop}\label{psat}
With the above notations, we have a natural Morita equivalence of $A$-linear dg-categories
$$\mathsf{MF}(X,\pi)^o \otimes_B \mathsf{MF}(X,\pi) \simeq L_{sg}(X\times_{S}X),$$
where the right hand side denotes the quotient $L_{\mathrm{coh}}^b(X\times_S X)/L_{\mathrm{perf}}(X\times_S X)$
computed in the Morita theory of dg-categories. 
\end{prop}

In the equivalence of Proposition \ref{psat}, the diagonal bimodule of $\mathsf{MF}(X,\pi)$ 
corresponds to the structure sheaf of the diagonal $\Delta_X \in L_{sg}(X\times_S X)$. 
This easily implies that $\mathsf{MF}(X,\pi)$ is a smooth and proper dg-category over
$B$. Moreover, the arguments developed in \cite{brtv}, together with Proposition 
\ref{psat}, imply that $\mathsf{MF}(X,\pi)$ satisfies the $r_{\ell}^{\otimes}$-admissibility 
conditions when the inertia group $I \subset Gal(K^{sp}/K)$ acts unipotently
on $H^*_{et}(X_{\overline{K}},\Ql)$. \\

\noindent Let us now come to our strategy for a possible proof of Bloch's Conjecture \ref{cb}. It 
is divided in two steps, and already the first one yields an interesting, and apparently new formula 
for the Artin conductor. \\

\noindent \textbf{First step: a conductor formula.} 
We first treat the case of unipotent monodromy.\\

\noindent \textsf{The unipotent case.}
We start by assuming the extra condition that the inertia subgroup $I \subset G_K=Gal(K^{sp}/K)$
acts with unipotent monodromy on all the $\Ql$-spaces $H^i_{\acute{e}t}(X_{\bar{K}},\Ql)$. This
implies that the monodromy action is in particular \emph{tame}, so that 
the Swan conductor term vanishes. Thus, Bloch's formula now reads
$$[\Delta_X,\Delta_X]_S = \chi(X_{\bar{k}},\ell) - \chi(X_{\bar{K}},\ell).$$

\begin{cor}\label{c1}
If the inertiia group $I$ acts unipotently, we have an equality
\begin{align*}
\chi(\mathsf{HH}(T/A)\otimes_{\mathsf{HH}(B/A)}B) = \chi(X_{\bar{k}},\ell) - \chi(X_{\bar{K}},\ell). \tag{CF-uni}
\end{align*}
\end{cor}

\noindent \textbf{Proof.} This is a direct consequence of our trace formula (Theorem \ref{t2}), and
the main theorem of \cite{brtv} that proves the existence of a natural equivalence of $\Ql(\beta)\oplus \Ql(\beta)[-1](1)$-modules
$$r_{\ell}(T) \simeq \mathbf{H}_{\acute{e}t}(X_{\bar{k}},\Phi_{X,\ell}[-1](\beta))^{hI},$$
where $\Phi_{X,\ell}$ is the complex of vanishing cycles on $X_k$, and $(-)^{hI}$ denotes the homotopy invariants for the action of the inertia group $I$. \hfill $\Box$ \\

\begin{rmk} (1) Let us remark that the $B$-linear dg-category $T=\mathsf{MF}(X,\pi)$ is not expected be $H^{\otimes}$-admissible
unless the monodromy action is unipotent. This is directly related to the fact that 
$r_{\ell}(T)$ provides $hI$-invariant vanishing cohomology, and taking invariants in general
does not commute with tensor products. However, when $I$ acts unipotently, taking 
$hI$-invariants does commute with tensor products computed over the dg-algebra $\Ql^{hI}=\Ql\oplus \Ql[-1](1)$. 
This is why Corollary \ref{c1} cannot be true for non-unipotent monodromy. \\
\noindent (2) Note that Bloch's conductor conjecture was open, 
without conditions on the reduced special fiber, 
even for unipotent monodromy. Therefore, Corollary \ref{c1} already provides
new cases where the conjecture is true.
\end{rmk}

\noindent \textsf{Extension to the non-unipotent case.} By Grothendieck's theorem (\cite[Exp. I]{sga7}), the action of the inertia $I$
is always quasi-unipotent. Let $S' \longrightarrow S$ be the totally ramified covering
corresponding to a totally ramified finite Galois extension $K'/K$ with group $G$, such that 
the base change $X':=X\times_S S' \longrightarrow S'$ has 
unipotent monodromy. The finite group $G$ acts on $X' \longrightarrow S'$ and induces
a morphism of quotient stacks
$$X'_G:=[X'/G] \longrightarrow S'_G:=[S'/G].$$
We now consider the category $T'$ of matrix factorizations
of $X'_G$ relative to $\pi$. This category is now linear over 
a $G$-equivariant $\E_2$-algebra $B_G$. The underlying $\E_2$-algebra 
$B_G$ is now $\mathbb{R}\underline{Hom}_{k\otimes_{A'} k}(k,k)[u^{-1}]$
where $S'=\Spec\, A'$. 

As the monodromy action for $X' \longrightarrow S'$ is unipotent, we conjecture that 
$T'$ is $H^{\otimes}$-admissible as a $G$-equivariant $B_G$-linear dg-category, and that it is smooth and proper, as well. This leads to a $G$-equivariant version
of Corollary \ref{c1} giving now a formula for the character of the $G$-representation
on $H^{*}_{\acute{e}t}(X',\nu[-1])^{hI}(\beta)$ in terms of the $G$-action on $\mathsf{HH}(T'/B_G)$. 
We hope that this $G$-equivariant version of Corollary \ref{c1} will imply the equality
\begin{align*}
\chi(\mathsf{HH}(\mathsf{MF}(X,\pi)/B)\otimes_{\mathsf{HH}(B/A)}B) = \chi(X_{\bar{k}}) - \chi(X_{\bar{K}}) - Sw(X_{\bar{K}}). \tag{CF-gen}
\end{align*}

\noindent \textbf{Second step: our conductor formula coincides with Bloch's conductor formula.} 
To complete the proof of Bloch's conductor Conjecture, we need to compare the localized intersection number
of Conjecture \ref{cb} and the Euler characteristic of the $B$-dg-module $\mathsf{HH}(T/B) \otimes_{\mathsf{HH}(B/A)}B$, where $T= \mathsf{MF}(X, \pi)$. This comparison should indeed be an equality
\begin{align*}
[\Delta_X,\Delta_X]_S = \chi(\mathsf{HH}(T/B) \otimes_{\mathsf{HH}(B/A)}B). \tag{Comp}
\end{align*}

\begin{rmk} When $A$ is equicharacteristic, then the twisted de Rham complex of $X$ (twisted by $\pi$) is defined (see e.g. \cite[Thm. 8.2.6]{pre}),  and can be used to prove the comparison (Comp). However, we think there is a different way to get formula (Comp) that works even in the mixed characteristic case, where the twisted de Rham complex is not  defined.
\end{rmk}

\bigskip
\bigskip

\noindent
Bertrand To\"{e}n, {\sc  Universit\'e Paul Sabati\'er \& CNRS},
Bertrand.Toen@math.univ-toulouse.fr 

\smallskip

\noindent
Gabriele Vezzosi, {\sc DIMAI, Universit\`a di Firenze},
gabriele.vezzosi@unifi.it

\end{document}